\documentclass{article}

\usepackage{amsmath, amsthm, amssymb}

\newtheorem{theorem}{Theorem}[section]
\newtheorem{lemma}[theorem]{Lemma}
\newtheorem{corollary}[theorem]{Corollary}
\theoremstyle{definition}
\newtheorem{definition}[theorem]{Definition}
\newtheorem{example}[theorem]{Example}

\newtheorem{problem}[theorem]{Problem}
\theoremstyle{remark}

\numberwithin{equation}{section}

\usepackage{graphicx}
\begin{document}

\title{Weyl-Heisenberg Frame Wavelets with Basic Supports}

\author{Xunxiang Guo, Yuanan Diao and Xingde Dai \\
Department of Mathematics }

%\address{Department of Mathematics}

\maketitle

\begin{abstract}
Let $a$, $b$ be two fixed non-zero constants. A measurable set
$E\subset \mathbb{R}$ is called a Weyl-Heisenberg frame set for
$(a, b)$ if the function $g=\chi_{E}$ generates a Weyl-Heisenberg
frame for $L^2(\mathbb{R})$ under modulates by $b$ and translates
by $a$, i.e., $\{e^{imbt}g(t-na\}_{m,n\in\mathbb{Z}}$ is a frame
for $L^2(\mathbb{R})$. It is an open question on how to
characterize all frame sets for a given pair $(a,b)$ in general.
In the case that $a=2\pi$ and $b=1$, a result due to Casazza and
Kalton shows that the condition that the set
$F=\bigcup_{j=1}^{k}([0,2\pi)+2n_{j}\pi)$ (where
$\{n_{1}<n_{2}<\cdots<n_{k}\}$ are integers) is a Weyl-Heisenberg
frame set for $(2\pi,1)$ is equivalent to the condition that the
polynomial $f(z)=\sum_{j=1}^{k}z^{n_{j}}$ does not have any unit
roots in the complex plane. In this paper, we show that this
result can be generalized to a class of more general measurable
sets (called basic support sets) and to set theoretical functions
and continuous functions defined on such sets.
\end{abstract}

\section{Introduction}

\medskip
Let $\mathbb{H}$ be a separable complex Hilbert space. Let
$B(\mathbb{H})$ denote the algebra of all bounded linear operators
on $\mathbb{H}$. Let $\mathbb{N}$ denote the set of natural
numbers, and $\mathbb{Z}$ be the set of all integers. A collection
of elements$ \{x_{j}:j\in \jmath\}$ in $\mathbb{H}$ is called a
\textit{frame} (of $\mathbb{H}$) if there exist constants $A$ and
$B$, $0<A\leq B< \infty$, such that
\begin{equation}
 A\|f\|^{2}\leq \sum_{j\in \jmath}|\langle f,\
x_{j}\rangle|^{2}\leq B\|f\|^{2}
\end{equation}
for all $f\in \mathbb{H}$. The supremum of all such numbers $A$
and the infimum of all such numbers $B$ are called the
\textit{frame bounds} of the frame and are denoted by $A_{0}$ and
$B_{0}$ respectively. A frame is called a \textit{tight frame} if
$A_{0}=B_{0}$ and is called a \textit{normalized tight frame} if
$A_{0}=B_{0}=1$. Any orthonormal basis in a Hilbert space is a
normalized tight frame. However a normalized tight frame is not
necessarily an orthonormal basis. Frames can be regarded as the
generalizations of orthogonal bases of Hilbert spaces. The concept
of frames was introduced a long time ago (\cite{DS,SZN}) and have
received much attention recently due to the development and study
of wavelet theory \cite{BMM,D,DB2,YM}. Among those widely studied
lately are the Weyl-Heisenberg frames (also called Gabor frames).
Let $a$, $b$ be two fixed positive constants and let $T_a$, $M_b$
be the \textit{translation operator by $a$} and \textit{modulation
operator by $b$} respectively, i.e., $T_af(t)=f(t-a)$ and
$M_bf(t)=e^{ibt}f(t)$ for any $f\in L^2(\mathbb{R})$. For a fixed
$g\in L^2(\mathbb{R})$, we say that $(g,a,b)$ generates a
\textit{Weyl-Heisenberg frame} if $\{M_{mb}T_{na}g\}_{m,n\in
\mathbb{Z}}$ forms a frame for $L^2(\mathbb{R})$. We also say that
the function $g$ is a \textit{mother Weyl-Heisenberg frame
wavelet} for $(a,b)$ in this case. Furthermore, a  measurable set
$E\subset \mathbb{R}$ is called a \textit{Weyl-Heisenberg frame
set} for $(a,\ b)$ if the function $g=\chi_E$ generates a
Weyl-Heisenberg frame for $L^2(\mathbb{R})$ under modulates by $b$
and translates by $a$, i.e., $\{e^{i mbt}g(t-na)\}_{m,n\in \mathbb
Z}$ is a frame for $L^2(\mathbb{R})$. Characterizing the mother
Weyl-Heisenberg frame wavelets in general is a difficult and open
question. In fact, even in the special case of $a=2\pi$ and $b=1$,
this question remains unsolved. There are many works related to
this subject, for more information please refer to
\cite{BC,C12,CC,CL,CCJ,CDH,FJ,SZ}.

\medskip
In 2002, Casazza and Kalton proved the following theorem \cite{C}:

\medskip
\begin{theorem}\label{t}
For fixed integers $n_1<n_2<\cdots <n_k$, the set
\begin{equation}\label{E1}
E=\cup_{j=1}^k\big([0,2\pi)+2\pi n_j\big)
\end{equation}
is a Weyl-Heisenberg frame set for $(2\pi,1)$ if and only if the
polynomial $p(z)=\sum_{j=1}^k z^{n_j}$ has no unit roots.
\end{theorem}

\medskip
This means that characterizing a Weyl-Heisenberg frame set $E$
with the special form (\ref{E1}) is equivalent to the following
problem, which was proposed by Littlewood in 1968 \cite{LITTLE}:

\begin{problem}\label{W} Classify the integer
sets $\{n_1<n_2<\cdots< n_k\}$ such that the polynomial
$p(z)=\sum_{1\le j\le k}z^{n_j}$ does not have any unit roots.
\end{problem}

While it is unfortunate that Theorem \ref{t} fails to provide a
definite answer to the characterization problem of Weyl-Heisenberg
frame sets since Problem \ref{W} is also an open question (see
\cite{BOW1,BOW2,C,LITTLE} and the references therein), it does
reveal a deep connection between the two seemingly irrelevant
subjects. In this paper, we will further investigate
Weyl-Heisenberg frames defined by set theoretical functions and
continuous functions with restricted domain. In Section 2, we will
introduce some concepts, definitions, preliminary lemmas. In
Section 3, we will state and prove our main theorems. Some
examples are given in Section 4.

\medskip
\section{Definitions and Preliminary Lemmas}

\begin{definition}\label{zakd}
The \textit{Zak transform} of a function $f\in L^2(\mathbb{R})$ is
defined as
\begin{equation}
Zf(t,w)=\frac{1}{\sqrt{2\pi}}\sum_{n\in \mathbb{Z}} f(t+2\pi n)e^{
i nw},\ \forall\ t,\ w \in [0,2\pi).
\end{equation}
\end{definition}

The above definition of Zak transform is slightly different from
that given in \cite{AJ}. We have the following lemma.

\begin{lemma}\cite{AJ}
The Zak transform is a unitary map from $L^2(\mathbb{R})$ onto
$L^2(Q)$, where $Q=[0, 2\pi)\times [0,2\pi)$.
\end{lemma}

One can easily check that the Zak transform maps

\begin{lemma}{\label{L1}}
\begin{equation}
Z(M_{m}T_{2n\pi}g)(t,w)=e^{i(mt+nw)}Zg(t,w),\ \forall\ g\in
L^2(\mathbb{R}).
\end{equation}
\end{lemma}
\proof Since $M_{m}T_{2n\pi}g(x)=e^{imx}g(x-2n\pi)$,
\begin{eqnarray*}
&&Z(M_{m}T_{2n\pi}g)(t,w)\\
&=&\frac{1}{\sqrt{2\pi}} \sum_{\ell \in \mathbb{Z}}e^{im(t+2\ell
\pi)}g\big(t+2(\ell-n)\pi\big)e^{i\ell w}\\
&=&\frac{1}{\sqrt{2\pi}} \sum_{\ell\in
\mathbb{Z}}e^{imt}g\big(t+2(\ell-n)\pi\big)e^{i\ell w}\\
&=&\frac{1}{\sqrt{2\pi}} \sum_{\ell\in
\mathbb{Z}}e^{i(mt+nw)}g\big(t+2(\ell-n)\pi\big)e^{i\ell
w-inw}\\
&=&\frac{1}{\sqrt{2\pi}}e^{i(mt+nw)}\sum_{\ell\in\mathbb{Z}}g\big(t+2(\ell-n)\pi\big)e^{i(\ell
-n)w}\\
&=&\frac{1}{\sqrt{2\pi}} e^{i(mt+nw)}\sum_{k\in
\mathbb{Z}}g(t+2k\pi)e^{ikw}\\
&=&e^{i(mt+nw)}\frac{1}{\sqrt{2\pi}}\sum_{k\in \mathbb{Z}}g(t+2k\pi)e^{ikw}\\
&=& e^{i(mt+nw)}(Zg)(t,w).
\end{eqnarray*}
\qed

\begin{lemma}\label{L2.4}\cite{C}
Let $E$ be a measurable subset of\ $[0,2\pi)$, $F=\bigcup_{n\in
\mathbb{Z}}(E+2n\pi)$ and $g\in L^2(F)$. The following statements
are equivalent:

\noindent (1) $(M_{m}T_{2n\pi}g)_{m,n\in \mathbb{Z}}$ is a frame
for $L^2(F)$ with frame bounds $A_0$, $B_0$.

\noindent (2) $0<A_0=essinf_{(t,w)\in E\times[0,2\pi)}
|Zg(t,w)|^2\leq esssup_{(t,w)\in E\times[0,2\pi)}
|Zg(t,w)|^2=B_0<\infty.$
\end{lemma}

Let $E=\bigcup_{k=1}^{m}A_{k}$, where each $A_{k}=[a_k,\ b_k]$ is
a closed interval and $A_i\cap A_j=\emptyset$ if $i\not=j$. Let
$\tau_{2\pi}(x):\ \mathbb{R}\longrightarrow [0,2\pi)$ be the
function defined by $\tau_{2\pi}(x)=x-2\pi[\frac{x}{2\pi}]$, where
$[.]$ is the integer function. Let us arrange the numbers
\[
0,\ \tau_{2\pi}(a_1),\ \tau_{2\pi}(b_1),\  \cdots,\
\tau_{2\pi}(a_m),\ \tau_{2\pi}(b_m),\ 2\pi
\]
in the ascending order and write the resulting numbers as $0=t_0<
t_1<\cdots< t_{j_{0}}=2\pi$. Then for each $0\le j\le j_0-1$ and
each $k$, we have either $(t_j,\ t_{j+1})\subset \tau_{2\pi}(A_k)$
or $(t_j,\ t_{j+1})\bigcap \tau_{2\pi}(A_k)=\emptyset$. Based on
this observation we have:
\begin{lemma}{\label{l2}}
Let $\{A_{i}\}_{i=1}^{m}$ be a sequence of finite and
non-overlapping intervals and $E=\bigcup_{i=1}^{m}A_{i}$, then
there exists a finite sequence of disjoint intervals
$\{E_{i}\}_{i=1}^{k}$ with $E_i\subset[0,2\pi)$, and an integer
sequence $\{n_{ij}\}_{j=1}^{j_{i}}$ for each $i$, such that
 \[
 E=\bigcup_{i=1}^{k}F_{i},\
\rm{where}\ \ F_{i}=\bigcup_{j=1}^{j_i}(E_{i}+2\pi n_{ij}).
\]
\end{lemma}

We will call the set $E$ defined in the above lemma a
\textit{basic  support set}, which is just a finite disjoint union
of finite intervals. The sequence $\{E_{i}\}_{i=1}^{k}$ associated
with the set $E$ will be called the \textit{2$\pi$-translation
generators } of $E$. Notice that $\bigcup_{j=1}^{j_i}(E_{i}+2\pi
n_{ij})$ is simply the pre-image of the function $\tau_{2\pi}$
restricted to $E$. We will call the sequence
$\{n_{i{j}}\}_{j=1}^{j_{i}}$ the \textit{step-widths} of the
corresponding generator $E_{i}$.

\begin{definition}
Let $g\in L^2(\mathbb{R})$ be a continuous function and $E$ be a
basic support set. For each $\xi\in E$, the sequence $\{g(\xi+2\pi
n_{i{j}})\}_{j=1}^{j_{i}}$ is called a \textit{characteristic
chain} of the function $g$ associated with the set $E$, where
$\xi\in E_i$ and $E_i$, $\{n_{ij}\}$ are as defined in Lemma
\ref{l2}.
\end{definition}

If $\sum_{j=0}^{k} a_{n_j}z^{n_j}=0$ has zeros on the unit circle
$\mathbb{T}$, then the coefficient sequence
$\{a_{n_j}\}_{j=0}^{k}$ is called a \textit{root sequence} of
$\{n_j\}$. For instance, $\{-2,-1,1\}$ is a root sequence for
$\{0,1,2\}$ since the polynomial $-2-z+z^2$ has a unit root. But
$\{-2,-1,1\}$ is not a root sequence for $\{0,1,3\}$ since the
polynomial $-2-z+z^3$ does not have a unit root. Furthermore, the
sequence $\{n_j\}$ may contain negative integers as well, in which
case $\sum_{j=0}^{k} a_{n_j}z^{n_j}$ is a Laurent polynomial, but
our results will still hold for such a polynomial.

\section{Main Results and their Proofs}

We first outline the main results obtained in this paper below.
The first theorem generalizes the result of \cite{C} to step like
functions.

\begin{theorem}{\label{r1}}
Let $n_{0}<n_1<n_2<\cdots<n_k$ be $k+1$ fixed integers, $a_{0}$,
$a_{1}$, $a_{2}$, $\cdots$, $a_{k}$ be $k+1$ given complex numbers
and let $F_{j}=[0,2\pi)+2\pi n_j$ ($j=0$, $1$, $2$, $\cdots$,
$k$), $\mbox{$g=\sum_{j=0}^{k}a_{j}\chi_{F_{j}}$}$. The following
statements are equivalent:

\smallskip
(1) $g$ is a mother Weyl-Heisenberg frame wavelet for $(2\pi,1)$
with frame bounds $A_0$, $B_0$.

\smallskip
(2) ${2\pi}A_0=\inf_{|z|=1}|\sum_{j=0}^{k}a_{j}z^{n_{j}}|^2$ and
${2\pi}B_0=\sup_{|z|=1}|\sum_{j=0}^{k}a_{j}z^{n_{j}}|^2$.

\smallskip
(3) For every measurable set $E\subset[0,2\pi)$ of positive
measure, let $E_{j}=E+2n_{j}\pi$, $j=0$, $1$, $2$, ..., $k$, the
function $g_E=\sum_{j=0}^{k}a_{j}\chi_{E_{j}}$ is a mother
Weyl-Heisenberg frame wavelet of $L^{2}(\Omega)$ for $(2\pi,1)$
with frame bounds $A\ge A_0$, $B\le B_0$, where
$\Omega=\cup_{n\in\mathbb{Z}}(E+2n\pi)$ and $A_0$, $B=B_0$ are the
frame bounds for $g_E$ with $E=[0,2\pi)$.
\end{theorem}

\begin{theorem}{\label{r3}}
Let $E\subset [0,2\pi)$ be a measurable set of positive measure
and $n_{0}<n_1<n_2<\cdots<n_k$ be $k+1$ given integers. Let
$F=\bigcup_{j=0}^{k}(E+2\pi n_{j})$, $\Omega =\bigcup_{n\in
\mathbb{Z}}(E+2\pi n)$. Then for any continuous function $g \in
L^2(\mathbb{R})$, $(g\cdot \chi_F, 2\pi,1)$ generates a
Weyl-Heisenberg frame of $L^2(\Omega)$ if and only if for any $\xi
\in \overline{E}$, $\{g(\xi+2\pi n_{j})\}$ is not a root sequence
for $\{n_j\}$ where $\overline{E}$ is the closure of $E$. In
particular, if $E=[0,2\pi)$, then $(g\cdot \chi_{F}, 2\pi,1)$
generates a Weyl-Heisenberg frame of $L^2(\mathbb{R})$ if and only
if for any $\xi \in [0,2\pi]$, $\{g(\xi+2\pi n_{j})\}$ is not a
root sequence for $\{n_j\}$.
\end{theorem}

Theorem \ref{r3} further generalizes Theorem \ref{r1} to
continuous functions.

\begin{theorem}{\label{r4}}
Let $E$ be a basic support set with generator and step-width pairs
$\{\big(E_{i},\ \{n_{i{j}}\}_{j=0}^{j_{i}} \big)\}_{i=0}^{k}$,
$\Omega=\bigcup_{n\in \mathbb{Z}}(E+2n\pi)$ and $g\in
L^2(\mathbb{R})$ be a continuous function, then $(g\cdot
\chi_E,2\pi,1)$ generates a Weyl-Heisenberg frame for
$L^2(\Omega)$ with frame bounds $A_0$ and $B_0$ if and only if
\[
2\pi A_0=\min_{0\le i\le k}\inf_{\xi\in \overline{E},z\in
\mathbb{T}} |\sum_{j=0}^{j_{i}}g(\xi+2\pi n_{i{j}})z^{n_{i{j}}}|^2
\]
and
\[
2\pi B_0=\max_{0\le i\le k}\sup_{\xi\in \overline{E},z\in
\mathbb{T}} |\sum_{j=0}^{j_{i}}g(\xi+2\pi
n_{i{j}})z^{n_{i{j}}}|^2.
\]
In other word, $g_E=g\cdot \chi_E$ is a mother Weyl-Heisenberg
frame wavelet of $L^2(\Omega)$ for $(2\pi,1)$ if and only if no
{characteristic chains} of $g$ associated with $\overline{E}$ are
root sequences (of the corresponding step-width sequences).
\end{theorem}

Notice that the special case $g=1$ of Theorem \ref{r4} will relate
the characterization of a Weyl-Heisenberg frame set within the
basic support sets (which is more general than the sets considered
in \cite{C}) to the classification of corresponding polynomials
with unit roots. The result in the following corollary to Theorem
\ref{r4} is immediate, which can also be obtained from the
definition of mother Weyl-Heisenberg frames directly without
difficulty.

\begin{corollary}
If $g\in L^2(\mathbb{R})$ is a continuous function with a zero
characteristic chain associated with a closed basic support set
$E$, then $g\cdot \chi_E$ is not a mother Weyl-Heisenberg frame
wavelet of $L^2(\Omega)$ for $(2\pi,1)$, where
$\Omega=\bigcup_{n\in \mathbb{Z}}(E+2n\pi)$.
\end{corollary}

\medskip
\textit{Proof of Theorem \ref{r1}.} For all $x,\ y\in[0,2\pi)$, we
have
\begin{eqnarray*}
& &Zg(x,y)=\frac{1}{\sqrt{2\pi}}\sum_{n\in \mathbb{Z}}g(x+2n\pi)e^{iny}\\
&=&\frac{1}{\sqrt{2\pi}}
\sum_{n\in\mathbb{Z}}\sum_{j=0}^{k}a_{j}\chi_{
F_{j}}(x+2n\pi)e^{iny}\\
&=&\frac{1}{\sqrt{2\pi}}
\chi_{[0,2\pi)}(x)\sum_{j=0}^{k}a_{j}e^{in_{j}y}.
\end{eqnarray*}

Thus
\begin{eqnarray*}
& &{\rm ess}\inf_{(x,y)\in[0,2\pi)^{2}}
|Zg(x,y)|^{2}\\
&=&\frac{1}{{2\pi}}{\rm ess}\inf_{y\in
[0,2\pi)}|\sum_{j=0}^{k}a_{j}e^{in_{j}y}|^{2}\\
&=&\frac{1}{{2\pi}}{\rm ess}\inf_{z\in
\mathbb{T}}|\sum_{j=0}^{k}a_{j}z^{n_{j}}|^{2}
\end{eqnarray*}
and
\begin{eqnarray*}
&&{\rm ess}\sup_{(x,y)\in[0,2\pi)^{2}}
|Zg(x,y)|^{2}\\
&=&\frac{1}{{2\pi}}{\rm ess}\sup_{y\in
[0,2\pi)}|\sum_{j=0}^{k}a_{j}e^{in_{j}y}|^{2}\\
& =&\frac{1}{{2\pi}}{\rm ess}\sup_{z\in
\mathbb{T}}|\sum_{j=0}^{k}a_{j}z^{n_{j}}|^{2},
\end{eqnarray*}
where $\mathbb{T}$ is the unit circle. By Lemma \ref{L2.4}, it
follows that $A_0$ and $B_0$ are the frame bounds of $g$ if and
only if ${2\pi}A_0=\inf_{|z|=1}|\sum_{j=0}^{k}a_{j}z^{n_{j}}|^2$
and ${2\pi}B_0=\sup_{|z|=1}|\sum_{j=0}^{k}a_{j}z^{n_{j}}|^2$. This
proves the equivalence of (1) and (2). The special case of
$E=[0,2\pi)$ in (3) implies (1). By (1), $(g,2\pi,1)$ generates a
Weyl-Heisenberg frame for $L^{2}(\mathbb{R})$ with frame bounds
$A$, $B$. Let $P$ be the orthogonal projection of
\[
L^{2}(\mathbb{R}) \rightarrow L^{2}(\Omega)
\]
defined by
\[
Pf=f|_{\Omega}=f\cdot\chi_{\Omega}.
\]
Then
\begin{eqnarray*}
&&P(E_{m}T_{2n\pi}g)=P\big(e^{imt}g(t-2n\pi)\big)\\
&=&\big(e^{imt}g(t-2n\pi)\big)\cdot\chi_{\Omega}\\
&=&e^{imt}g_E(t-2n\pi)=E_{m}T_{2n\pi}g_E.
\end{eqnarray*}
So $g_E=\sum_{j=0}^{k}a_{j}\chi_{ E_{j}}$ is a mother
Weyl-Heisenberg frame for $L^{2}(\Omega)$ with frame bounds $A\ge
A_0$, $B\le B_0$.

\qed

\noindent{\it  Proof of Theorem \ref{r3}.} Since $g$ is continuous
on $\mathbb{R}$, $g_F=g\cdot \chi_F$ is continuous on $F$ and it
follows that $Zg_F(x,y)$ is continuous on $E\times [0,2\pi)$. So
there exist $0<A\le B$ such that
\[
A={\rm ess}\inf_{(x,y)\in E\times [0,2\pi)}|Zg(x,y)|^2\leq {\rm
ess}\sup_{(x,y)\in E\times [0,2\pi)}|Zg(x,y)|^2= B
\]
 if and only if
\[
0<A=\inf_{x\in E,\ z\in \mathbb{T}}
|\frac{1}{\sqrt{2\pi}}\sum_{j=0}^{k}g(x+2\pi n_{j})z^{n_{j}}|^2
\]
and
\[
B=\sup_{x\in E,\ z\in \mathbb{T}}
|\frac{1}{\sqrt{2\pi}}\sum_{j=0}^{k}g(x+2\pi n_{j})z^{n_{j}}|^2
\]
by Lemma \ref{L2.4}. The result of the Theorem then follows since
$A$ and $B$ are actually attained on the set $\overline{E}\times
\mathbb{T}$ when $g$ is continuous. \qed

\medskip
\noindent{\it Proof of Theorem \ref{r4}.} By the definition of
basic support sets, we have $E=\bigcup_{i=1}^{k}F_{i}$, where
\[
F_{i}=\bigcup_{j=1}^{j_{i}}(E_{i}+2n_{i{j}}\pi)
\]
and the $E_i$'s are disjoint intervals in $[0,2\pi)$. Let
$M_{i}=\bigcup_{n\in \mathbb{Z}}(E_{i}+2n\pi)$, then $
M_{i}\bigcap M_{j}=\phi \ for\ any\ i\ne j $,
\[
\Omega=\bigcup_{n\in \mathbb{Z}}(E+2n\pi)=\bigcup_{i=1}^{k}M_{i}
\]
and
\[
L^2(\Omega)=\bigoplus_{i=1}^{k}L^2(M_{i}).
\]
It is easy to see that $g$ is a mother Weyl-Heisenberg frame
wavelet of $L^2(\Omega)$ if and only if for each $i$, $g\cdot
\chi_{F_{i}}$ is a mother Weyl-Heisenberg frame wavelet of
$L^2(M_{i})$. By Theorem \ref{r3}, $g\cdot \chi_{F_{i}}$ is a
mother Weyl-Heisenberg frame wavelet of $L^2(M_{i})$ if and only
if no characteristic chains of $g$ associated with
$\overline{F_i}$ are root sequences (of their corresponding
step-width sequences). We leave it to our reader to check that
\[
A_0=\frac{1}{2\pi}\min_{0\le i\le k}\inf_{\xi\in \overline{E},z\in
\mathbb{T}} |\sum_{j=0}^{j_{i}}g(\xi+2\pi n_{i{j}})z^{n_{i{j}}}|^2
\]
and
\[
B_0=\frac{1}{2\pi}\max_{0\le i\le k}\sup_{\xi\in \overline{E},z\in
\mathbb{T}} |\sum_{j=0}^{j_{i}}g(\xi+2\pi n_{i{j}})z^{n_{i{j}}}|^2
\]
are the frame bounds of $g\cdot \chi_E$ in case that $A_0>0$
(hence $g\cdot \chi_E$ is a mother Weyl-Heisenberg frame wavelet
of $L^2(\Omega)$). \qed

\bigskip
\section{Examples}

\medskip
\begin{example}
The set $[3\pi,7\pi)$ is a Weyl-Heisenberg frame set: We have
$E_{1}=[0,\pi)$, $E_2=[\pi,2\pi)$, $n_{1{1}}= 2$, $n_{1{2}}=3$,
$n_{2{1}}=1$, $n_{2{2}}=2$ since
\[[3\pi,7\pi)=\big(E_1+4\pi\big)\cup\big(E_1+6\pi\big)\cup\big(E_2+2\pi\big)\cup\big(E_2+4\pi\big).
\]
Furthermore, the corresponding polynomials $2+3z$ and $1+2z$ have
no unit roots.
\end{example}

\begin{example}
The set
\[
(\frac{5}{2}\pi,\frac{7}{2}\pi]\cup(4\pi,\frac{11}{2}\pi]=\big((0,\frac{1}{2}\pi]+4\pi\big)\cup\big((\frac{1}{2}\pi,
\frac{3}{2}\pi]+2\pi\big)\cup\big((\frac{1}{2}\pi,
\frac{3}{2}\pi]+4\pi\big)
\]
is a WH-frame set since $E_{1}=(0,\frac{1}{2}\pi]$,
$E_2=(\frac{1}{2}\pi, \frac{3}{2}\pi]$, $n_{1{1}}=2$,
$n_{2{1}}=1$, $n_{2{2}}=2$ and the polynomials $2$ and $1+2z$ have
no unit roots.
\end{example}

\begin{example}
On the unit circle of complex plane,
\[
|4+3z+2z^3|=|2+3z^2+4z^3|.
\]
The roots of $p(z)=4+3z+2z^3$ are: $r_{1}=0.4398+1.4423i$,
$r_{2}=0.4398-1.4423i$, $r_{3}=-0.8796$. So $p(z)$ doesn't have
unit zeros on the complex plane. It follows that the set functions
\[
g_{1}(\xi)=4\chi_{[0,2\pi)}+3\chi_{[2\pi,4\pi)}+2\chi_{[6\pi,8\pi)}
\]
and
\[
g_{2}(\xi)=2\chi_{[0,2\pi)}+3\chi_{[4\pi,6\pi)}+4\chi_{[6\pi,8\pi)}
\]
are both mother Weyl-Heisenberg frame wavelets for $(2\pi,1)$ with
the same frame bounds.
\end{example}

\begin{example}
Similarly, using the fact that $
|1-2z+3z^3-5z^5|=|-5+3z^2-2z^4+z^5|$ on the unit circle and that
$p(z)=1-2z+3z^3-5z^5$ doesn't have unit zeros, we see that the
following set functions
\[
g_{3}(\xi)=\chi_{[0,2\pi)}-2\chi_{[2\pi,4\pi)}+3\chi_{[6\pi,8\pi)}-5\chi_{[10\pi,12\pi)}
\]
and
\[
g_{4}(\xi)=-5\chi_{[0,2\pi)}+3\chi_{[4\pi,6\pi)}-2\chi_{[8\pi,10\pi)}+\chi_{[10\pi,12\pi)}
\]
are both mother Weyl-Heisenberg frame wavelets for $(2\pi,1)$ with
the same frame bounds.
\end{example}

\begin{example}
The function $f(t)=\sin t$ is not a mother Weyl-Heisenberg frame
wavelet of $L^2(\mathbb{R})$ for $(2\pi,1)$ on any basic support
set $E$ (in fact any measurable set) since for any set $E$ such
that $\mathbb{R}=\cup_{n\in \mathbb{Z}}(E+2n\pi)$, $\overline{E}$
must contain a subsequence of $\{k\pi\}_{k\in \mathbb{Z}}$.
Consequently, the function $\sin t$ will always have a zero
characteristic chain associated with the set $\overline{E}$.
However, if $E$ is a WH-frame set such that $\overline{E}$ is
disjoint from $\{k\pi\}_{k\in \mathbb{Z}}$, then $\sin t\cdot
\chi_E$ is a mother WH-frame wavelet of $L^2(\Omega)$ for
$(2\pi,1)$ where $\Omega=\cup_{n\in \mathbb{Z}}(E+2n\pi)$ since
any characteristic chain of $\sin t$ associated with
$\overline{E}$ is a constant sequence. This example can be
generalized to other $2\pi$-periodic functions. In particular,
given any continuous $2\pi$-periodic function $g(t)$ that is
bounded away from $0$ and a WH-frame set $E$, $g\cdot \chi_E$ is a
mother WH-frame wavelet of $L^2(\mathbb{R})$ for $(2\pi,1)$. For
instance, the function
\[
g(t)=\left\{\begin{array}{ll} |\sin(t)| & {\rm if}\ t\in
[\frac{\pi}{6},\frac{5\pi}{6})+k\pi,\ \forall k\in \mathbb{Z},\\
\frac{1}{2} & {\rm otherwise.}
\end{array}
\right.
\]
is a mother Weyl-Heisenberg frame wavelet of $L^2(\mathbb{R})$ for
$(2\pi,1)$ when it is restricted to the set $E=[3\pi,7\pi)$. See
Figure \ref{fig1} below.

\begin{figure}[!htb]
\begin{center}
\includegraphics[angle=270,scale=0.4]{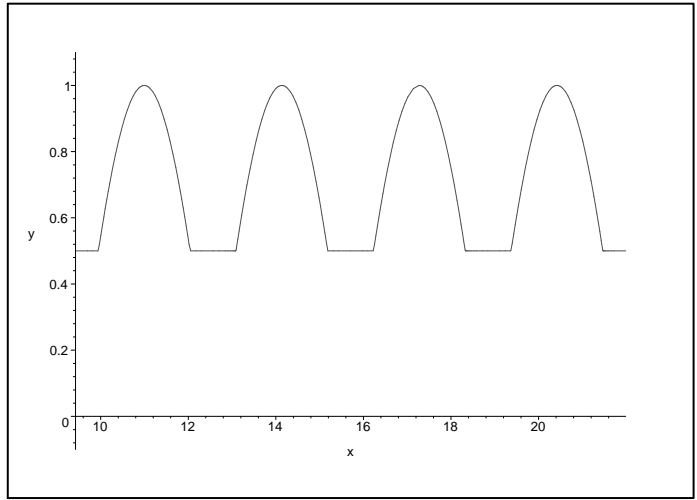}
\caption{\label{fig1} A continuous periodic function bounded away
from $0$ on a WH-frame set.}
\end{center}
\end{figure}
\end{example}

\begin{example}
Show that the following function $g(t)$ (see Figure \ref{fig2}) is
a mother Weyl-Heisenberg frame wavelet of $L^2(\mathbb{R})$ for
$(2\pi,1)$ when it is restricted to the set $E=[0,2\pi)\cup
[4\pi,6\pi)\cup [8\pi,10\pi)$. Notice that the set $E$ is not a
WH-frame set.
\begin{figure}[!htb]
\begin{center}
\includegraphics[angle=270,scale=0.5]{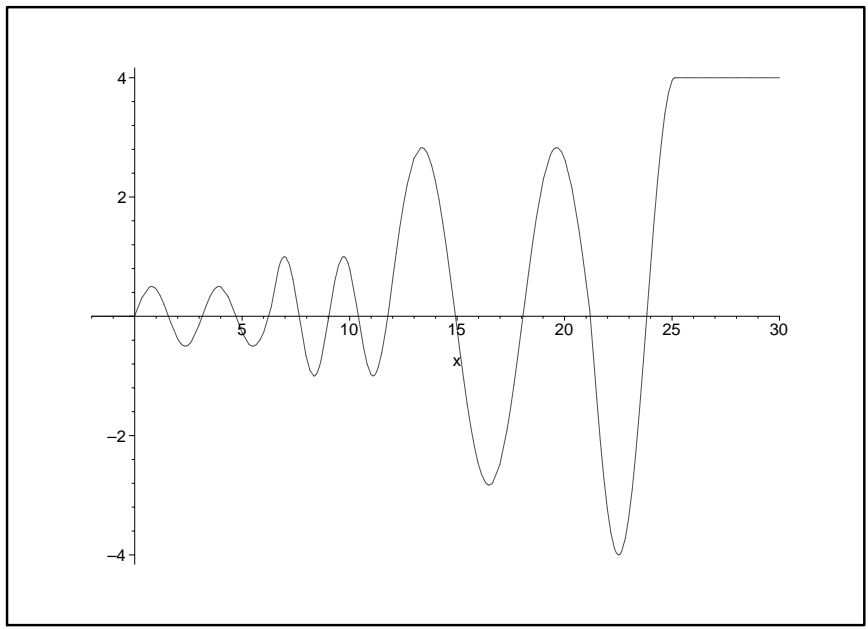}
\caption{\label{fig2} A continuous non-periodic function that is
not bounded away from $0$ which defines a mother Weyl-Heisenberg
frame wavelet of $L^2(\mathbb{R})$ for $(2\pi,1)$ when it is
restricted to a set that is not a WH-frame set.}
\end{center}
\end{figure}
The function $g(t)$ is defined piecewisely by
\[
g(t)=\left\{\begin{array}{ll}
0 & t<0,\\
\sin(2t)/2 & t\in [0,2\pi),\\
\sin(\frac{16}{7}(t-2\pi)) & t\in [2\pi, \frac{15}{4}\pi),\\
2(\sin t+\cos t),& t\in [\frac{15}{4}\pi, \frac{27}{4}\pi),\\
-4\sin(\frac{6}{5}(t-\frac{27}{4}\pi)),& t\in
[\frac{27}{4}\pi,8\pi),\\
 4 & t\ge 8\pi.
\end{array}
\right.
\]
In this case, $E_1=[0,2\pi)$ and $n_{11}=0$, $n_{12}=2$,
$n_{13}=4$ since $E=E_1\cup(E_1+2\cdot 2\pi)\cup (E_1+4\cdot
2\pi)$. So for any given $t\in E_1$, the corresponding
characteristic chain of $g(t)$ with respect to $E$ is
$\{\sin(2t)/2,2(\sin t + \cos t), 4\}$. It is not a root sequence
of $\{0,2,4\}$ since $\sin(2t)/2+2(\sin t+\cos t)z^2+4z^4=(\sin
t+2z^2)(\cos t+2z^2)$ apparently have no unit zeros for any given
$t$.
\end{example}

\bigskip

\bibliographystyle{amsalpha}

\end{document}